\documentclass[12pt,twoside]{article}
\usepackage{times}
\usepackage{amsmath,amssymb}
\usepackage{color}
\usepackage{fancyhdr}
\allowdisplaybreaks

 \def \nn{\nonumber}

\newcommand{\pf}{\noindent {\bf Proof. \hspace{2mm}}}
\newcommand{\ef}{ \hfill $ \blacksquare $ \vskip 3mm}
\newcommand{\be}{\begin{equation}}
\newcommand{\ee}{\end{equation}}
\newcommand{\bea}{\begin{eqnarray}}
\newcommand{\eea}{\end{eqnarray}}

\newcommand{\bR}{{\mathbb R}}

\def\p{\partial}
\def\la{\lambda}

\def\al{\alpha}

\def\q{\quad}

\def\g{\gamma}

\def\dl{\delta}
\def\ve{\varepsilon}

\def\lt{\left}
\def\rt{\right}

\def\q{\qquad}

\def\dl{\delta}

\def\p{\partial}
\def\f{\frac}

\def\al{\alpha}

\def\q{\qquad}
\def\s{\sqrt}
 \allowdisplaybreaks

\begin{document}
\let\oldsection\section
\renewcommand\section{\setcounter{equation}{0}\oldsection}
\renewcommand\thesection{\arabic{section}}
\renewcommand\theequation{\thesection.\arabic{equation}}
\newtheorem{claim}{\noindent Claim}[section]
\newtheorem{theorem}{\noindent Theorem}[section]
\newtheorem{lemma}{\noindent Lemma}[section]
\newtheorem{proposition}{\noindent Proposition}[section]
\newtheorem{definition}{\noindent Definition}[section]
\newtheorem{remark}{\noindent Remark}[section]
\newtheorem{corollary}{\noindent Corollary}[section]
\newtheorem{example}{\noindent Example}[section]

\title{Remarks on 1-D Euler Equations with Time-Decayed Damping}

\author{Xinghong Pan
\footnote{E-mail:math.scrat@gmail.com.}\vspace{0.5cm}\\
\small  (Department of Mathematics and
IMS, Nanjing University, Nanjing 210093, China.)\\
\vspace{0.5cm}
}

\date{}
\maketitle

\centerline {\bf Abstract} \vskip 0.3 true cm

We study the 1-d isentropic Euler equations with time-decayed damping
\be
\left\{
\begin{aligned}
&\p_t \rho+\p_x(\rho u)=0, \\
&\p_t(\rho u)+ \p_x(\rho u^2)+\p_xp(\rho)=-\f{\mu}{1+t}\rho u,\\
&\rho|_{t=0}=1+\ve\rho_0(x),u|_{t=0}=\ve u_0(x).
\end{aligned}
\right. \nn
\ee

 This work  is inspired by a recent work of F. Hou, I. Witt and H.C. Yin \cite{Hou01}. In \cite{Hou01}, they proved a global existence and blow-up result of 3-d  irrotational Euler flow with time-dependent damping. In the 1-d case, we will prove a different result when the damping decays of order $-1$ with respect to the time $t$. More precisely, when $\mu>2$, we prove the global existence of the 1-d Euler system. While when $0\leq\mu\leq2 $, we will prove the blow up of $C^1$ solutions.

\vskip 0.3 true cm

\vskip 0.3 true cm


\section{Introduction}
\q This paper deals with the isentropic Euler equations with time-decayed damping in 1 dimension:
\be
\left\{
\begin{aligned}
&\p_t \rho+\p_x(\rho u)=0, \\
&\p_t(\rho u)+ \p_x(\rho u^2)+\p_xp(\rho)=-\f{\mu}{1+t}\rho u,\\
&\rho|_{t=0}=1+\ve\rho_0(x),u|_{t=0}=\ve u_0(x),
\end{aligned}
\right.\label{e1.1}
\ee
where $\rho_0(x), u_0(x)\in C^\infty_0(\bR)$, supported in $|x|\leq R$ and $\ve>0$ is sufficient small. Here $\rho(x),u(x)$ and $p(x)$ represent the density, fluid velocity and pressure respectively and $\mu$ is a positive number to describe the scale of the damping. We assume the fluid is a polytropic gas which means we assume $p(\rho)=\f{1}{\g}\rho^\g, \g>1$. We denote $c^2=p'(\rho).$

As is well known, when the damping is vanishing, the smooth solution of compressible Euler flow will blow up in finite time. For the extensive literature on the blow-up results and the blow-up mechanism, Readers can see \cite{D01}, \cite{D02}, \cite{T01,T02}, and  \cite{Y01,Y02} and references therein for more details.

While 1-d Euler equations with linear damping read as
\be
\left\{
\begin{aligned}
&\p_t \rho+\p_x(\rho u)=0, \\
&\p_t(\rho u)+ \p_x(\rho u^2)+\p_xP(\rho)=-\kappa\rho u.\\
\end{aligned}
\right.  \label{1.2}
\ee

Here $\kappa$ is a constant. There are many results to prove the global existence and convergence rates of solutions to system \eqref{1.2} with small data. Readers can turn to\cite{02}, \cite{09} and \cite{10} for more information.

In multi space dimensions, there are also some results to Euler equations with linear damping. Wang and Yang proved the stability of the planar diffusion wave for the 2D Euler system with linear damping in \cite{14} and studied the pointwise estimates of solutions for the 3-dimension Euler equations in \cite{15}. Similar convergence rates in 3 dimensions are obtained by Tan and Wang by using a different method in
\cite{12}. Also in \cite{04}, Jiu and Zheng proved the global existence of the 3D Euler systems with linear damping in Besov spaces. While Sideris, Thomases and Wang in \cite{11} showed that the smooth solutions of the linear-damped Euler equations do not decay in exponentially in time and may blow up if the initial data is sufficient large.

It is nature to ask whether the global solution exists when the damping is decayed and what is the critical decayed rate to separate the global existence and the finite-time blow up of solutions with small data. Recently F. Hou, I. Witt and H.C. Yin in \cite{Hou01} proved a global existence(the damping is time-decayed smaller than or equal to order -1) and blow-up(the damping is time-decayed larger than order -1) result for the 3d irrotational Euler flow. Our article discusses the 1-dimension Euler equations with damping of time-decayed order -1. Compared with their work, we will get a different result in 1-dimensional case. The global existence and blow up depend on the scale of $\mu$ and $\mu=2$ is the critical value.

Throughout this paper we denote a generic constant by $C$. It may be different line by line. $H^m(R)$ denotes the usual Sobolev space with its norm
\be
\|f\|_{H^m}\triangleq \sum\limits^m_{k=0}\|\p^k_xf\|_{L_p}. \nn
\ee
For convenience, we use $\|\cdot\|$ to denote $\|\cdot\|_{L_2}$ and $\|\cdot\|_m$ for $\|\cdot\|_{H^m}$.

We state our main results as following.
\begin{theorem}
(Global existence for $ \mu>2$) Suppose $(\rho_0,u_0)\in H^m(\bR), m\geq 3$ and $\mu>2$. Then there exists a unique global classical solution $(\rho(x,t),u(x,t))$ of \eqref{e1.1} satisfies
\bea
&&\|(1+t)\p_t\rho(t)\|^2_{m-1}+\|(1+t)\p_x\rho(t)\|^2_{m-1}+\|(1+t)\p_xu(t)\|^2_{m-1} \nn  \\
&&+\|(\rho(t)-1)\|^2 +\|u(t)\|^2    \nn  \\
&&+\int^t_0(1+\tau)(\|\p_\tau\rho(\tau)\|^2_{m-1}+\|\p_x\rho(\tau)\|_{m-1}+\|\p_x u(\tau)\|_{m-1})d\tau    \nn  \\
&\leq& C\ve^2(\|\rho_0\|^2_{m}+\|u_0\|^2_{m}).  \nn
\eea
\end{theorem}
\begin{remark}
By Sobolev inequality, we have the pointwise estimates
\be
\sum\limits_{\al\leq m-1}\sup\lt\{(1+t)^{\f{1+\min\{\al,1\}}{2}}(|\p^\al_xv|+|\p^\al_xu|)\rt\}\leq C\ve. \nn
\ee
\end{remark}
\begin{remark}
Inspired by the work of F. Hou, I. Witt and H.C. Yin, we believe that by introducing a ``potential function'' $\phi(x,t)=\int^x_\infty u(y,t)dy$, we can prove the global existence when the damping  decays like $\f{\mu}{(1+t)^\la}$ $(0<\la<1, \mu>0)$ . We will present it in our later work.
\end{remark}

The idea of proving Theorem1.1 comes from \cite{15}. In the first part of their paper, they proved the global existence of multi-dimensional Euler equations with constant linear damping by a energy estimate method. We revise their method with a time-weighted energy estimate to make it suitable for the proving of Theorem1.1.

Next we discuss the blow up of system \eqref{e1.1} when  $ 0\leq\mu\leq 2.$ Define $T_\ve$ be the lifespan, the largest existing time, of $C^1$ solutions to system \eqref{e1.1}.
Define two functions
\be
q^0(r)=\int^\infty_r(x-r)^2\rho_0(x)dx, \nn
\ee
\be
q^1(r)=2\int^\infty_r(x-r)(\rho_0u_0)(x)dx. \nn
\ee
\begin{theorem}
Suppose for some $0<R_0<R,$
\be
q^0(r)>0,\q q^1(r)\geq 0,  \label{1.3}
\ee
for $R_0<r<R.$ Then the lifespan $T_\ve$ of $C^1$ solutions of \eqref{e1.1} is finite.
\end{theorem}
\begin{remark}
the proof of Theorem1.2 can find its original version in T. Sideris \cite{T02} for 3-d compressible Euler equations. It can not be used directly for the 1-d damped Euler equations, but we can revise the test function there to similarly prove the finite-time blow up.
\end{remark}
\begin{remark}
Actually our method can deal with the case when the damping decays like $\f{\mu}{(1+t)^\la}$ when $\la>1, \mu\geq 0.$
\end{remark}

We arrange our paper as following. In Section 2, we prove Theorem1.1 for the global existence with relatively ``large'' damping. In Section 3, we prove Theorem 1.3 for the blow up of solutions with relatively ``small'' damping.
\section{PROOF OF THEOREM 1.1}

\q In this section, we prove Theorem 1.1 by a time-weighted energy method. Remember $c=\s{P'(\rho)}=\rho^{\f{\g-1}{2}}$. First we change \eqref{e1.1} into
\be
\left\{
\begin{aligned}
&\f{2}{\g-1}\p_t c+c\p_x u+\f{2}{\g-1}u \p_xc=0, \\
&\p_t u+ u\p_x u+\f{2}{\g-1}c\p_xc+\f{\mu}{1+t} u=0,\\
&c|_{t=0}=1+\ve c_0(x),u|_{t=0}=\ve u_0(x),
\end{aligned}
\right.\label{e2.1}
\ee
where $c_0(x)\in C^\infty_0(R)$, supported in $|x|\leq R.$

Let $v=\f{2}{\g-1}(c-1)$, then $v, u$ satisfy
\be
\left\{
\begin{aligned}
&\p_t v+\p_x u=-u \p_xv-\f{\g-1}{2}v\p_xu, \\
&\p_t u+ \p_xv+\f{\mu}{1+t} u=-u\p_xu-\f{\g-1}{2}v\p_xv,\\
&v|_{t=0}=\ve v_0(x),u|_{t=0}=\ve u_0(x),
\end{aligned}
\right.\label{e2.2}
\ee
where $v_0(x)=\f{2}{\g-1}c_0(x).$

From \eqref{e2.2}, we have

\be
\p_{tt} v-\p_{xx}v+\f{\mu}{1+t}\p_tv=Q(v,u),\label{e2.3}
\ee
where
\bea
Q(v,u)&=&\f{\mu}{1+t}(-u\p_xv-\f{\g-1}{2}v\p_xu)             \nn            \\
&-&\p_t(u\p_xv)-\f{\g-1}{2}\p_t(v\p_xu)+\p_x(u\p_xu)+\f{\g-1}{2}\p_x(v\p_xv).  \nn
\eea
In the following, we will estimate $(v,u)$ under the a priori assumption
\bea
E_m(T)&=&:\sup\limits_{0<t<T}\lt\{\|(1+t)\p_xv\|^2_{m-1}+\|(1+t)\p_xu\|^2_{m-1}+\|v\|^2+\|u\|^2\rt\}^{\f{1}{2}} \nn \\
     &\leq& M\ve.         \label{e2.4}
\eea
Where $M$, independent of $\ve$, will be determined later. By choosing $M$ suitably, we will prove
\be
E_m(T)\leq \f{1}{2}M\ve.  \label{aa}
\ee
By Sobolev inequality, we know that
\bea
\sum\limits_{|\al|\leq 1}&&\sup\lt\{(1+t)^{\f{\al+1}{2}}(|\p^\al_xv|+|\p^\al_x u|)\rt\}   \nn  \\
&&\leq CE_2(t)\leq CM\ve.   \nn
\eea
Since \eqref{e2.2} implies
\be
|(1+t)\p_t v|+|(1+t)\p_t u|\leq C\lt\{|u|+(1+t)|\p_xv|+(1+t)|\p_xu|\rt\}, \nn
\ee
we also have
\be
\sum\limits_{\al+\beta\leq 1}\sup\lt\{(1+t)^{\f{\al+\beta+1}{2}}(|\p^\al_x\p^\beta_t v|+|\p^\al_x\p^\beta_t u|)\rt\}\leq CM\ve. \label{e2.5}
\ee

In the following, we will first obtain some elementary estimates for the 1-order derivatives of the solution. Then the higher derivatives will be handled in the similar way.

\subsection{Estimate 1}

For some constant $\eta$, to be determined later, multiplying \eqref{e2.3} by $\eta(1+t)^2\p_t v$ yields
\bea
&&\f{\eta}{2}\p_t[(1+t)^2(\p_tv)^2]+\eta(\mu-1)(1+t)(\p_tv)^2-\eta(1+t)^2\p_tv\p_{xx}v  \nn  \\
&=&\eta(1+t)^2\p_tvQ(v,u).     \label{e2.6}
\eea
Integrating it over $R\times[0,t]$ and using integration by parts give
\bea
&&\f{\eta}{2}\int_R(1+t)^2(\p_tv)^2dx+\eta(\mu-1)\int^t_0\int_R(1+\tau)(\p_\tau v)^2dxd\tau  \nn  \\
&&+\f{\eta}{2}\int_R(1+t)^2(\p_xv)^2dx-\eta\int^t_0\int_R(1+\tau)(\p_xv)^2dxd\tau   \nn  \\
&=&\f{\eta}{2}\int_R(\p_tv)^2|_{t=0}dx+\f{\eta}{2}\int_R(\p_xv)^2|_{t=0}dx+\int^t_0\int_R\eta(1+\tau)^2\p_\tau vQ(v,u)dxd\tau  \nn  \\
&\leq& C\ve^2(\|v_0\|^2_{1}+\|u_0\|^2_{1})+\int^t_0\int_R\eta(1+\tau)^2\p_\tau vQ(v,u)dxd\tau.   \label{e2.7}
\eea
Also multiplying \eqref{e2.3} by $(1+t)v$, we get
\bea
&&\p_t[(1+t)v\p_tv]+\f{\mu-1}{2}\p_tv^2-(1+t)v\p_{xx}v-(1+t)(\p_tv)^2  \nn  \\
&=&(1+t)vQ(v,u).   \label{e2.8}
\eea
Then integrating \eqref{e2.8} over $R\times[0,t]$ and using integration by parts give
\bea
&&\int_R(1+t)v\p_tvdx+\f{\mu-1}{2}\int_Rv^2dx      \nn  \\
&&+\int^t_0\int_R(1+\tau)(\p_xv)^2dxd\tau-\int^t_0\int_R(1+\tau)(\p_\tau v)^2dxd\tau   \nn  \\
&=&\int_R(v\p_tv)|_{t=0}dx+\f{\mu-1}{2}\int_Rv^2|_{t=0}dx+\int^t_0\int_R(1+t)vQ(v,u)dxd\tau  \nn  \\
&\leq& C\ve^2(\|v_0\|^2_{1}+\|u_0\|^2_{1})+\int^t_0\int_R(1+\tau)vQ(v,u)dxd\tau.   \label{e2.9}
\eea
Adding \eqref{e2.7} and \eqref{e2.9}, we have
\bea
&&\f{\eta}{2}\int_R(1+t)^2(\p_tv)^2dx+\f{\eta}{2}\int_R(1+t)^2(\p_xv)^2dx+\f{\mu-1}{2}\int_Rv^2dx  \nn \\
&&+[\eta(\mu-1)-1]\int^t_0\int_R(1+\tau)(\p_\tau v)^2dxd\tau+(1-\eta)\int^t_0\int_R(1+\tau)(\p_xv)^2dxd\tau  \nn  \\
&&+\int_R(1+t)v\p_tvdx   \nn \\
&\leq& C\ve^2(\|v_0\|^2_{1}+\|u_0\|^2_{1})+\int^t_0\int_R[(1+\tau)^2\p_\tau v+(1+\tau)v]Q(v,u)dxd\tau.   \label{e2.10}
\eea

If $\mu=2+4\dl$ for some $\dl>0$, using Cauchy-Schwartz inequality, we have
\bea
\int_R(1+t)v\p_tvdx\geq -\f{1+4\dl}{2(1+\dl)}\int_Rv^2dx-\f{1+\dl}{2(1+4\dl)}\int_R(1+t)^2(\p_tv)^2dx.  \label{e2.11}
\eea
From \eqref{e2.10} and \eqref{e2.11} by choosing $\eta=\f{1+2\dl}{1+4\dl}$, we have
\bea
&&\f{\dl}{2(1+4\dl)}\int_R(1+t)^2(\p_tv)^2dx+\f{1+2\dl}{2(1+4\dl)}\int_R(1+t)^2(\p_xv)^2dx+\f{\dl(1+4\dl)}{2(1+\dl)}\int_Rv^2dx  \nn \\
&&+2\dl\int^t_0\int_R(1+\tau)(\p_\tau v)^2dxd\tau+\f{2\dl}{1+4\dl}\int^t_0\int_R(1+\tau)(\p_xv)^2dxd\tau  \nn  \\
&\leq& CE^2_1(0)+\int^t_0\int_R[(1+\tau)^2\p_\tau v+(1+\tau)v]Q(v,u)dxd\tau.   \label{e2.12}
\eea
So we have
\bea
&&\int_R(1+t)^2(\p_tv)^2dx+\int_R(1+t)^2(\p_xv)^2dx+\int_Rv^2dx  \nn \\
&&+\int^t_0\int_R(1+\tau)(\p_tv)^2dxd\tau+\int^t_0\int_R(1+\tau)(\p_xv)^2dxd\tau  \nn  \\
&\leq& CE^2_1(0)+CI.   \label{e2.13}
\eea
where $C$ depends on $\mu$ and
\bea
I&=&\int^t_0\int_R[(1+\tau)^2\p_\tau v+(1+\tau)v]Q(v,u)dxd\tau  \nn  \\
 &=&\int^t_0\int_R[(1+\tau)^2\p_\tau v+(1+\tau)v]\f{\mu}{1+\tau}(-u\p_xv-\f{\g-1}{2}v\p_xu)dxd\tau     \nn                   \\
 &&-\int^t_0\int_R[(1+\tau)^2\p_\tau v+(1+\tau)v][\p_\tau(u\p_xv)+\f{\g-1}{2}\p_\tau(v\p_xu)]dxd\tau   \nn   \\
&&+\int^t_0\int_R[(1+\tau)^2\p_\tau v+(1+\tau)v][\p_x(u\p_xu)+\f{\g-1}{2}\p_x(v\p_xv)]dxd\tau   \nn  \\
 &=&I_1+I_2+I_3.       \nn
\eea
Now we estimate $I_1, I_2$ and $I_3$.
\bea
I_1&=&\mu\int^t_0\int_R(1+\tau)\p_\tau v(-u\p_xv-\f{\g-1}{2}v\p_xu)dxd\tau  \nn    \\
    &&+\mu\int^t_0\int_Rv(-u\p_xv-\f{\g-1}{2}v\p_xu)dxd\tau.   \nn
\eea
From \eqref{e2.5}, using Cauchy-Schwartz inequality and integration by parts, we have
\bea
I_1&\leq& CM\ve \int^t_0\int_R\lt[(1+\tau)(|\p_\tau v|^2+|\p_xv|^2+|\p_xu|^2)+\f{u^2}{1+\tau}\rt]dxd\tau  \nn  \\
   &&+C\int^t_0\int_R uv\p_x vdxd\tau    \nn  \\
   &\leq&CM\ve \int^t_0\int_R\lt[(1+\tau)(|\p_\tau v|^2+|\p_xv|^2+|\p_xu|^2)+\f{u^2}{1+\tau}\rt]dxd\tau. \label{e2.14}
\eea

Now we focus on the estimate of $I_2$, then $I_3$ will be essentially the same with $I_2$.
Dealing $I_2$ the same with $I_1$, we have
\bea
I_2&\leq& CM\ve\int^t_0\int_R(1+\tau)(|\p_\tau u|^2+|\p_xv|^2+|\p_\tau v|^2)dxd\tau   \nn   \\
&&-\int^t_0\int_R(1+\tau)^2u\p_\tau v\p^2_{x\tau}vdxd\tau-\int^t_0\int_R(1+\tau)vu\p^2_{x\tau}vdxd\tau    \nn  \\
&&-\f{\g-1}{2}\int^t_0\int_R(1+\tau)^2v\p_\tau v\p^2_{x\tau}udxd\tau-\f{\g-1}{2}\int^t_0\int_R(1+\tau)v^2\p^2_{x\tau}udxd\tau  \nn \\
&=&CM\ve\int^t_0\int_R(1+\tau)(|\p_\tau u|^2+|\p_xv|^2+|\p_\tau v|^2)dxd\tau   \nn   \\
&&+\f{1}{2}\int^t_0\int_R(1+\tau)^2\p_xu(\p_\tau v)^2dxd\tau+\int^t_0\int_R(1+\tau)\p_x(vu)\p_\tau vdxd\tau    \nn  \\
&&-\f{\g-1}{2}\int^t_0\int_R(1+\tau)^2v\p_\tau v\p^2_{x\tau}udxd\tau+\f{\g-1}{2}\int^t_0\int_R(1+\tau)\p_xv^2\p_\tau udxd\tau.  \nn
\eea
Using \eqref{e2.5} again, we have
\bea
I_2&\leq&CM\ve\int^t_0\int_R(1+\tau)(|\p_\tau u|^2+|\p_xv|^2+|\p_\tau v|^2)dxd\tau   \nn   \\
&&-\f{\g-1}{2}\int^t_0\int_R(1+\tau)^2v\p_\tau v\p^2_{x\tau}udxd\tau.  \label{e2.15}
\eea
From \eqref{e2.2}
\be
\p_xu=\f{-\p_tv-u\p_xv}{1+\f{\g-1}{2}v}.  \nn  \label{2.16}
\ee
Then
\be
\p^2_{xt}u=-\f{\p^2_tv+\p_tu\p_xv+u\p^2_{xt}v}{1+\f{\g-1}{2}v}+\f{\f{\g-1}{2}\p_tv(\p_tv+u\p_xv)}{(1+\f{\g-1}{2}v)^2}. \label{e2.16}
\ee
From \eqref{e2.5}, we have
\be
\f{1}{1+\f{\g-1}{2}v}\leq \f{1}{1-CM\ve}.    \label{e2.17}
\ee
Inserting \eqref{e2.16} and \eqref{e2.17} into \eqref{e2.15}, we have
\bea
&&\int^t_0\int_R(1+\tau)^2v\p_\tau v\p^2_{x\tau}udxd\tau  \nn   \\
&=&\int^t_0\int_R(1+\tau)^2v\p_\tau v\lt[-\f{\p^2_\tau v+\p_\tau u\p_xv+u\p^2_{x\tau}v}{1+\f{\g-1}{2}v}+\f{\f{\g-1}{2}\p_\tau v(\p_\tau v+u\p_xv)}{(1+\f{\g-1}{2}v)^2}\rt]dxd\tau \nn   \\
&\leq&-\int^t_0\int_R(1+\tau)^2\f{v}{1+\f{\g-1}{2}v}\p_\tau v(\p^2_\tau v+u\p^2_{x\tau}v)dxd\tau    \nn   \\
&&+\f{CM^2\ve^2}{(1-CM\ve)^2}\int^t_0\int_R (1+\tau)(|\p_\tau v|^2+|\p_xv|^2)dxd\tau.  \label{e2.18}
\eea
Combining \eqref{e2.15} and \eqref{e2.18}, we have
\bea
I_2&\leq&C\lt[M\ve+\f{M^2\ve^2}{(1-CM\ve)^2}\rt]\int^t_0\int_R (1+\tau)(|\p_\tau v|^2+|\p_xv|^2+|\p_\tau u|^2)dxd\tau  \nn \\
&&-\int^t_0\int_R\lt[(1+\tau)^2\f{v}{1+\f{\g-1}{2}v}\p_\tau v\p^2_\tau v+(1+\tau)^2\f{vu}{1+\f{\g-1}{2}v}\p_\tau v\p^2_{x\tau}v\rt]dxd\tau \nn  \\
&=&C\lt[M\ve+\f{M^2\ve^2}{(1-CM\ve)^2}\rt]\int^t_0\int_R (1+\tau)(|\p_\tau v|^2+|\p_xv|^2+|\p_\tau u|^2)dxd\tau  \nn \\
&&+I_{2,1}+I_{2,2}.  \label{e2.19}
\eea
Using integration by parts with respect to $\tau$, we have
\bea
I_{2,1}&=&-\f{1}{2}\int^t_0\int_R(1+\tau)^2\f{v}{1+\f{\g-1}{2}v}\p_\tau (\p_\tau v)^2dxd\tau  \nn  \\
       &=&-\f{1}{2}\int_R(1+t)^2\f{v}{1+\f{\g-1}{2}v} (\p_t v)^2dx+\f{1}{2}\int_R\f{\ve v_0}{1+\f{\g-1}{2}\ve v_0} (\p_\tau v(\cdot,0))^2dx \nn \\
       &&+\f{1}{2}\int^t_0\int_R\p_\tau\lt[(1+\tau)^2\f{v}{1+\f{\g-1}{2}v}\rt] (\p_\tau v)^2dxd\tau \nn   \\
       &\leq&\f{C\ve}{1-C\ve}E^2_1(0)+\f{CM\ve}{1-CM\ve}\int_R(1+t)^2(\p_tv)^2dx \nn  \\
       &&+\f{CM\ve}{(1-CM\ve)^2}\int^t_0\int_R(1+\tau) (\p_\tau v)^2dxd\tau. \label{e2.20}
\eea
And using integration by parts, we have
\bea
I_{2,2}&=&-\f{1}{2}\int^t_0\int_R(1+\tau)^2\f{vu}{1+\f{\g-1}{2}v}\p_x(\p_\tau v)^2dxd\tau  \nn  \\
       &=&\f{1}{2}\int^t_0\int_R(1+\tau)^2\p_x\lt[\f{vu}{1+\f{\g-1}{2}v}\rt](\p_\tau v)^2dxd\tau  \nn  \\
       &\leq&\f{CM^2\ve^2}{(1-CM\ve)^2}\int^t_0\int_R(1+\tau)^2(\p_\tau v)^2dxd\tau.  \label{e2.21}
\eea
Inserting \eqref{e2.20} and \eqref{e2.21} into \eqref{e2.19}, we have
\bea
I_2&\leq&\f{CM\ve}{(1-CM\ve)^2}\int^t_0\int_R (1+\tau)(|\p_\tau v|^2+|\p_xv|^2+|\p_\tau u|^2)dxd\tau  \nn  \\
   &&+\f{C\ve}{1-C\ve}E^2_1(0)+\f{CM\ve}{1-CM\ve}\int_R(1+t)^2(\p_tv)^2dx.  \label{e2.22}
\eea
We can deal $I_3$ almost the same with $I_2$. Then we can get
\bea
I_3&\leq&\f{CM\ve}{(1-CM\ve)^2}\int^t_0\int_R (1+\tau)(|\p_\tau v|^2+|\p_xv|^2+|\p_x u|^2)dxd\tau  \nn  \\
   &&+\f{C\ve}{1-C\ve}E_1(0)+\f{CM\ve}{1-CM\ve}\int_R(1+t)^2(\p_xv)^2dx.  \label{e2.23}
\eea
Remember that in \eqref{e2.2}, we have
\be
|\p_t u|\leq C\lt(|\p_x v|+|\p_x u|+\f{|u|}{1+t}\rt).  \nn
\ee
Inserting the estimates of $I_1$ \eqref{e2.14}, $I_2$ \eqref{e2.22} and $I_3$ \eqref{e2.23} into \eqref{e2.13}, we have
\bea
&&\int_R(1+t)^2(\p_tv)^2dx+\int_R(1+t)^2(\p_xv)^2dx+\int_Rv^2dx  \nn \\
&&+\int^t_0\int_R(1+\tau)(\p_\tau v)^2dxd\tau+\int^t_0\int_R(1+\tau)(\p_xv)^2dxd\tau  \nn  \\
&\leq& C\f{1+\ve}{1-C\ve}E^2_1(0)+ \f{CM\ve}{1-CM\ve}\int_R  (1+t)^2\lt[(\p_\tau v)^2+(\p_xv)^2\rt]dx \nn  \\
&&+ \f{CM\ve}{(1-CM\ve)^2}\int^t_0\int_R\lt\{(1+\tau)\lt[(\p_\tau v)^2+(\p_xv)^2+(\p_xu)^2\rt]+\f{u^2}{1+\tau}\rt\}dxd\tau. \nn \\ \label{e2.24}
\eea
\subsection{Estimate 2}

Multiplying $\eqref{e2.2}_2$ by $u$ and integrating on $R\times[0,t]$ yield
\bea
\int^t_0\int_R(\f{1}{2}\p_\tau u^2+u\p_xv+\f{\mu}{1+\tau}u^2)dxd\tau=\int^t_0\int_R(-u\p_xu-\f{\g-1}{2}v\p_xv)udxd\tau. \nn
\eea
Then using integration by parts, we have
\bea
&&\f{1}{2}\int_R u^2dx-\f{\ve^2}{2}\int_R u^2_0dx-\int^t_0\int_Rv\p_xudxd\tau+\int^t_0\int_R\f{\mu}{1+\tau}u^2dxd\tau \nn  \\
&=&\int^t_0\int_R(-u\p_xu-\f{\g-1}{2}v\p_xv)udxd\tau.  \nn
\eea

Using \eqref{e2.2} and \eqref{e2.5}, we have
\bea
&&\f{1}{2}\int_R u^2dx+\int^t_0\int_Rv(\p_\tau v+u\p_xv+\f{\g-1}{2}v\p_xu)dxd\tau+\int^t_0\int_R\f{\mu}{1+\tau}u^2dxd\tau  \nn \\
&\leq&CE^2_1(0)+CM\ve\int^t_0\int_R\lt\{\f{1}{1+\tau}u^2+(1+\tau)\lt[(\p_xv)^2+(\p_xu)^2\rt]\rt\}dxd\tau. \nn
\eea
Then we can have
\bea
&&\int_R(u^2+v^2)dx+\int^t_0\int_R\f{\mu}{1+\tau}u^2dxd\tau \nn   \\
&\leq& CE^2_1(0)+CM\ve\int^t_0\int_R\lt\{\f{1}{1+\tau}u^2+(1+\tau)\lt[(\p_xv)^2+(\p_xu)^2\rt]\rt\}dxd\tau. \label{e2.25}
\eea
\subsection{Estimate 3}
By differentiating $\eqref{e2.2}_2$ with respect to $x$ and integrating its product with $(1+\tau)^2\p_xu$ on $R\times[0,t]$, we have
\bea
&&\int^t_0\int_R(1+\tau)^2\p_xu(\p^2_{x\tau}u+\p^2_xv+\f{\mu}{1+\tau}\p_xu)dxd\tau \nn  \\
&=&\int^t_0\int_R(1+\tau)^2\p_xu\p_x(-u\p_xu-\f{\g-1}{2}v\p_xv)dxd\tau.  \nn
\eea
Then
\bea
&&\f{1}{2}\int_R(1+t)^2(\p_xu)^2dx+(\mu-1)\int^t_0\int_R(1+\tau)(\p_xu)^2dxd\tau  \nn  \\
&&+\int^t_0\int_R(1+\tau)^2\p_x(\p_\tau v+u\p_xv+\f{\g-1}{2}v\p_xu)\p_xvdxd\tau \nn  \\
&\leq& CE^2_1(0)+\int^t_0\int_R(1+\tau)^2\p_xu\p_x(-u\p_xu-\f{\g-1}{2}v\p_xv)dxd\tau.  \nn
\eea
So
\bea
&&\f{1}{2}\int_R(1+t)^2\lt[(\p_xu)^2+(\p_xv)^2\rt]dx  \nn  \\
&&+(\mu-1)\int^t_0\int_R(1+\tau)(\p_xu)^2dxd\tau-\int^t_0\int_R(1+\tau)(\p_xv)^2dxd\tau \nn  \\
&\leq&-\int^t_0\int_R(1+\tau)^2\p_x(u\p_xv+\f{\g-1}{2}v\p_xu)\p_xvdxd\tau  \nn  \\
&&+CE^2_1(0)+\int^t_0\int_R(1+\tau)^2\p_xu\p_x(-u\p_xu-\f{\g-1}{2}v\p_xv)dxd\tau.  \nn
\eea
We can deal with the right terms almost the same with $I_2$ and $I_3$. Then we can get the estimate
\bea
&&\int_R(1+t)^2\lt[(\p_xu)^2+(\p_xv)^2\rt]dx \nn  \\
&&+\int^t_0\int_R(1+\tau)(\p_xu)^2dxd\tau -\int^t_0\int_R(1+\tau)(\p_xv)^2dxd\tau \nn  \\
&\leq& C\f{1+\ve}{1-C\ve}E^2_1(0)+ \f{CM\ve}{1-CM\ve}\int_R  (1+t)^2\lt[(\p_tv)^2+(\p_xv)^2\rt]dx \nn  \\
&&+ \f{CM\ve}{(1-CM\ve)^2}\int^t_0\int_R\lt\{(1+\tau)\lt[(\p_\tau v)^2+(\p_xv)^2+(\p_xu)^2\rt]+\f{u^2}{1+\tau}\rt\}dxd\tau. \nn \\ \label{e2.26}
\eea

Let \eqref{e2.24}$+$\eqref{e2.25}$+\f{1}{2}\times$\eqref{e2.26}, then we get
\bea
&&\int_R\lt\{(1+t)^2\lt[(\p_tv)^2+(\p_xv)^2+(\p_xu)^2\rt]+v^2+u^2\rt\}dx \nn  \\
&&+\int^t_0\int_R\lt\{(1+\tau)\lt[(\p_\tau v)^2+(\p_xv)^2+(\p_xu)^2\rt]+\f{u^2}{1+\tau}\rt\}dxd\tau \nn \\
&\leq& C\f{1+\ve}{1-C\ve}E^2_1(0)+ \f{CM\ve}{1-CM\ve}\int_R  (1+t)^2\lt[(\p_tv)^2+(\p_xv)^2\rt]dx \nn  \\
&&+ \f{CM\ve}{(1-CM\ve)^2}\int^t_0\int_R\lt\{(1+\tau)\lt[(\p_\tau v)^2+(\p_xv)^2+(\p_xu)^2\rt]+\f{u^2}{1+\tau}\rt\}dxd\tau. \nn \\ \label{e2.27}
\eea

\subsection{Estimates for Higher Derivatives}

The estimates as \eqref{e2.24} and \eqref{e2.26} can also be obtained for higher derivatives. In fact, by multiplying \eqref{e2.3} by $\p^2_x[\eta(1+\tau)^2\p_\tau v+(1+\tau)v]$ and integrating it on $R\times[0,t]$, we can get
\bea
&&\int_R(1+t)^2(\p^2_{xt}v)^2dx+\int_R(1+t)^2(\p^2_xv)^2dx+\int_R(\p_xv)^2dx  \nn \\
&&+\int^t_0\int_R(1+\tau)(\p^2_{x\tau}v)^2dxd\tau+\int^t_0\int_R(1+\tau)(\p^2_xv)^2dxd\tau  \nn  \\
&\leq& C\f{1+\ve}{1-C\ve}E^2_2(0)+ \f{CM\ve}{1-CM\ve}\int_R  (1+t)^2\lt[(\p^2_{xt}v)^2+(\p^2_xv)^2\rt]dx \nn  \\
&&+ \f{CM\ve}{(1-CM\ve)^2}\int^t_0\int_R\lt\{(1+\tau)\lt[(\p^2_{x\tau}v)^2+(\p^2_xv)^2+(\p^2_xu)^2\rt]\rt\}dxd\tau. \nn \\ \label{e2.28}
\eea
By differentiating \eqref{e2.2}$_2$ two times with respect to $x$ and integrating its product with $(1+\tau)^2\p^2_xu$ on $R\times[0,t]$, we have
\bea
&&\int_R(1+t)^2\lt[(\p^2_xu)^2+(\p^2_xv)^2\rt]dx+\int^t_0\int_R(1+\tau)(\p^2_xu)^2dxd\tau  \nn  \\
&\leq& C\f{1+\ve}{1-C\ve}E^2_2(0)+ \f{CM\ve}{1-CM\ve}\int_R  (1+t)^2\lt[(\p^2_{xt}v)^2+(\p^2_xv)^2\rt]dx \nn  \\
&&+ \f{CM\ve}{(1-CM\ve)^2}\int^t_0\int_R\lt\{(1+\tau)\lt[(\p^2_{x\tau}v)^2+(\p^2_xv)^2+(\p^2_xu)^2\rt]\rt\}dxd\tau.\nn \\ \label{e2.29}
\eea
Combing \eqref{e2.27}, \eqref{e2.28} and \eqref{e2.29} gives
\bea
&&\|(1+t)\p_tv\|^2_1+\|(1+t)\p_xv\|^2_1+\|(1+t)\p_xu\|^2_1+\|v\|^2+\|u\|^2  \nn  \\
&&+\int^t_0\lt[(1+\tau)(\|\p_\tau v\|^2_1+\|\p_xv\|^2_1+\|\p_xu\|^2_1)+\f{\|u\|^2}{1+\tau}\rt]d\tau  \nn \\
&\leq&C\f{1+\ve}{1-C\ve}E^2_2(0)+ \f{CM\ve}{1-CM\ve}(\|(1+t)\p_tv\|^2_1+\|(1+t)\p_xv\|^2_1) \nn  \\
&&+ \f{CM\ve}{(1-CM\ve)^2}\int^t_0\lt[(1+\tau)(\|\p_\tau v\|^2_1+\|\p_xv\|^2_1+\|\p_xu\|^2_1)+\f{\|u\|^2}{1+\tau}\rt]d\tau.\nn \\
\eea
Actually, we can prove for $m$
\bea
&&\|(1+t)\p_tv\|^2_{m-1}+\|(1+t)\p_xv\|^2_{m-1}+\|(1+t)\p_xu\|^2_{m-1}+\|v\|^2+\|u\|^2  \nn  \\
&&+\int^t_0\lt[(1+\tau)(\|\p_\tau v\|^2_{m-1}+\|\p_xv\|^2_{m-1}+\|\p_xu\|^2_{m-1})+\f{\|u\|^2}{1+\tau}\rt]d\tau  \nn \\
&\leq&C\f{1+\ve}{1-C\ve}E^2_m(0)+ \f{CM\ve}{1-CM\ve}(\|(1+t)\p_tv\|^2_{m-1}+\|(1+t)\p_xv\|^2_{m-1}) \nn  \\
&&+ \f{CM\ve}{(1-CM\ve)^2}\int^t_0\lt[(1+\tau)(\|\p_\tau v\|^2_{m-1}+\|\p_xv\|^2_{m-1}+\|\p_xu\|^2_{m-1})+\f{\|u\|^2}{1+\tau}\rt]d\tau. \nn \\
\eea
When $\ve$ is small, for some $C_0$, we get
\bea
&&\|(1+t)\p_tv\|^2_{m-1}+\|(1+t)\p_xv\|^2_{m-1}+\|(1+t)\p_xu\|^2_{m-1}+\|v\|^2+\|u\|^2  \nn  \\
&&+\int^t_0\lt[(1+\tau)(\|\p_\tau v\|^2_{m-1}+\|\p_xv\|^2_1+\|\p_xu\|^2_{m-1})+\f{\|u\|^2}{1+\tau}\rt]d\tau  \nn \\
&\leq&\f{(1-CM\ve)^2}{(1-CM\ve)^2-CM\ve}C_0\ve^2. \nn
\eea
Let $M^2=5C_0$. By using the smallness of $\ve$, we can have
\bea
E^2_m(t)\leq \f{1}{4}M^2\ve^2.   \label{e2.33}
\eea

The local existence of symmetrizable hyperbolic equations have been proved by using the fixed point theorem. In order to get the global existence of the system, we only need a priori estimate. Based on our above estimate \eqref{e2.33} and the continuation argument, we finish the prove of Theorem 1.1.\ef
\section{PROOF OF THEOREM 1.2}

\q In this section, we prove Theorem1.2 when $0\leq\mu\leq 2$.
We first deal with the case $\g=2$ and later indicate the modification for general case.

\pf
Let $(\rho, u)$ be a $C^1$ solution. By the finite propagation property, we have $\rho-1$ supported in $B(t)=\{x\mid|x|\leq t+R\}.$ We define
\be
P(r,t)=\int_{x>r}(x-r)^2(\rho(x,t)-1)dx. \label{3.1}
\ee
Using \eqref{e1.1}$_1$ and integration by parts, we have
\bea
\p_t P(r,t)&=&\int_{x>r}(x-r)^2\rho_tdx \nn  \\
&=&-\int_{x>r}(x-r)^2(\rho u)_xdx  \nn  \\
&=&\int_{x>r}2(x-r)(\rho u)dx. \nn
\eea
Then $P(r,t)$ is $C^2$ in $t$. Differentiating it again, using \eqref{e1.1}$_2$ and integration by parts, we have
\bea
\p^2_t P(r,t)&=&\int_{x>r}2(x-r)(\rho u)_t dx  \nn  \\
&=&\int_{x>r}2(x-r)\lt(-\p_x(\rho u^2)-\p_x p-\f{\mu}{1+t}\rho u\rt) dx  \nn  \\
&=&\int_{x>r}2\rho u^2dx+\int_{x>r}2 pdx+\f{\mu}{1+t}\int_{x>r}(x-r)^2(\rho u)_x dx  \nn  \\
&=&\int_{x>r}2\rho u^2dx+\int_{x>r}2 pdx-\f{\mu}{1+t}\p_t\int_{x>r}(x-r)^2(\rho-1) dx.  \nn
\eea
Hence we have
\bea
\p^2_t P(r,t)+\f{\mu}{1+t}\p_t P(r,t)&=&\int_{x>r}2\rho u^2dx+\int_{x>r}2 pdx  \nn \\
&\geq 0.&          \nn
\eea

Due to our initial data assumption \eqref{1.3}, by integrating the above differential inequality, we have $\p_t P(r,t)\geq 0$ and $P(r,t)>0.$

Now we come to estimate a lower bound for $P(r,t)$. Rewriting $\p^2_t P(r,t)$ as following.
\bea
\p^2_t P(r,t)&=&\int_{x>r}2(x-r)(\rho u)_t dx  \nn  \\
&=&\int_{x>r}2(x-r)\lt(-\p_x(\rho u^2)-\p_x p-\f{\mu}{1+t}\rho u\rt) dx  \nn  \\
&=&\int_{x>r}2\rho u^2dx+\int_{x>r}2 (p-\overline{p})dx+\f{\mu}{1+t}\int_{x>r}(x-r)^2(\rho u)_x dx  \nn  \\
&\geq&\int_{x>r}2(p-\overline{p})dx-\f{\mu}{1+t}\p_t\int_{x>r}(x-r)^2(\rho-1) dx,  \nn
\eea
where $\overline{p}=p(1).$ Then we have
\bea
\p^2_t P(r,t)-\p^2_rP(r,t)+\f{\mu}{1+t}\p_tP(r,t)&\geq&\int_{x>r}2(p-\overline{p})dx-\int_{x>r}2(\rho-1)dx  \nn \\
&=&\f{2}{\g}\int_{x>r}\lt[ {(\rho^\g-1})-\g(\rho-1)\rt]dx  \nn  \\
&\triangleq&G(r,t). \label{e3.2}
\eea

When $\g=2,$
\be
G(r,t)=\int_{x>r}(\rho-1)^2dx\geq 0.  \nn
\ee

 When $0\leq \mu\leq 2,$ due to the nonnegativity of $\p_t P$, we have
\be
\p^2_t P(r,t)-\p^2_rP(r,t)+\f{2}{1+t}\p_tP(r,t)\geq G(r,t).\label {e3.2}\footnote{If the damping decays like $\f{\mu}{(1+t)^\la}$ where $\la>1,\mu\geq 0$, we can choose a $t_1$ such that when $t\geq t_1$, $\f{\mu}{(1+t)^\la}\leq \f{2}{1+t}$. We can still get \eqref{e3.2}.}
\ee
Define $W(r,t)=(1+t)P(r,t)$. From the above inequality, one get
\be
\p^2_t W(r,t)-\p^2_rW(r,t)\geq (1+t)G(r,t). \label{3.3}
\ee
We see
\be
W(r,0)=\ve q^0(r), \q (\p_tW)(r,0)=\ve(q^0(r)+q^1(r)).  \nn
\ee
Inversion of 1-d d'Alembertian operator gives (for $r>R_0+t$)
\bea
W(r,t)&=& W^0(r,t)+\f{1}{2}\int^t_0\int^{r+t-\tau}_{r-(t-\tau)}\Box W(y,\tau)dyd\tau  \nn  \\
&\geq& W^0(r,t)+\f{1}{2}\int^t_0\int^{r+t-\tau}_{r-(t-\tau)} (1+\tau)G(y,\tau)dyd\tau,
\eea
where
\be
W^0(r,t)=\f{\ve}{2}\lt\{q^0(r+t)+q^0(r-t)+\int^{r+t}_{r-t}(q^0(y)+q^1(y))dy\rt\}.  \label{3.5}
\ee
Now define
\be
F(t)=\int^t_0(t-\tau)\int^{R+\tau}_{R_0+\tau}r^{-1}W(r,t)drd\tau.  \nn
\ee
We see that
\bea
F^{\prime\prime}(t)&=&\int^{R+t}_{R_0+t}r^{-1}W(r,t)dr \nn  \\
 &\geq& \int^{R+t}_{R_0+t}r^{-1}W^0(r,t)dr  \nn   \\
 &&+\f{1}{2}\int^{R+t}_{R_0+t}r^{-1}\int^t_0\int^{r+t-\tau}_{r-(t-\tau)}(1+\tau)G(y,\tau)dyd\tau dr  \nn  \\
 &=& J_1+J_2.   \label{3.6}
\eea
From our assumption \eqref{1.3}, we have
\bea
J_1&\geq&\ve\f{1}{2}\int^{R+t}_{R_0+t}r^{-1}q^0(r-t)dr \nn  \\
   &\geq&\ve(R+t)^{-1}\f{1}{2}\int^{R+t}_{R_0+t}q^0(r-t)dr  \nn  \\
   &\geq&\ve(R+t)^{-1}B_0  \nn \\
   &>&0,  \label{3.7}
\eea
where $B_0=\f{1}{2}\int^R_{R_0}q^0(r)dr$.

Exchanging the order of integration in $J_2$ and remembering that $G(y,\tau)$ is supported in $\{y\mid |y|\leq \tau+R\},$ we have
\be
J_2\geq\f{1}{2}\int^t_0\int^{\tau+R}_{\tau+R_0}(1+\tau)G(y,\tau)\int^{y+t-\tau}_{\max[t+R_0,y-(t-\tau)]}r^{-1}drdyd\tau. \nn
\ee
If we set $t\geq t_1=\f{1}{2}(R-R_0),$ by direct computation, we have
\be
\int^{y+t-\tau}_{\max[t+R_0,y-(t-\tau)]}r^{-1}dr\geq C(t+R)^{-2}(t-\tau)(y-\tau-R_0)^2.
\ee

Since $G(y,\tau)\geq 0,$ we have
\bea
J_2\geq C(t+R)^{-2}\int^t_0\int^{\tau+R}_{\tau+R_0}(t-\tau)(y-\tau-R_0)^2(1+\tau)G(y,\tau)dyd\tau, \label{3.9}
\eea
when $t>t_1.$ We know that $G(y,\tau)$ is supported in $\{y\mid |y|\leq \tau+R\}$ and
\be
G(y,\tau)=\p^2_y\int_{x>y}(x-y)^2(\rho-1)^2dx.  \nn
\ee
Using integration by parts in \eqref{3.9}, we have
\bea
J_2&\geq& C(t+R)^{-2}\int^t_0(t-\tau)\int^{\tau+R}_{\tau+R_0}\int_{x>y}(1+\tau)(x-y)^2(\rho-1)^2dxdyd\tau \nn  \\
   &=& C(t+R)^{-2}J_3.  \label{3.10}
\eea
Recall that
\be
F(t)=\int^t_0(t-\tau)\int^{\tau+R}_{\tau+R_0}y^{-1}(1+\tau)\int_{x>y}(x-y)^2(\rho-1)dxdyd\tau.  \nn
\ee
Using Cauchy-Schwartz's inequality, we have
\bea
F^2(t)&\leq& J_3\int^t_0(t-\tau)(1+\tau)\int^{\tau+R}_{\tau+R_0}y^{-2}\int^{\tau+R}_y(x-y)^2dxdyd\tau   \nn  \\
      &=&J_3J_4.  \label{3.11}
\eea
We compute $J_4$ as follows
\bea
J_4&=&\f{1}{3}\int^t_0(t-\tau)(1+\tau)\int^{\tau+R}_{\tau+R_0}y^{-2}(\tau+R-y)^3dyd\tau  \nn \\
&\leq& \f{(R-R_0)^3}{3}\int^t_0(t-\tau)(1+\tau)\int^{\tau+R}_{\tau+R_0}y^{-2}dyd\tau  \nn  \\
&\leq& C\int^t_0(t-\tau)(1+\tau)\f{1}{(\tau+R_0)^2}d\tau  \nn  \\
&\leq& C(t+R)\ln(t+R). \label{3.12}
\eea
Combining \eqref{3.10}, \eqref{3.11}, \eqref{3.12} and \eqref{3.6}, we get
\be
F^{\prime\prime}(t)\geq C[(t+R)^3\ln(t+R)]^{-1}F^2(t), \q t\geq t_1.  \label{3.13}
\ee
From \eqref{3.6}, \eqref{3.7} and the fact $J_2\geq 0$, $F^\prime(0)=F(0)=0,$ we have
\bea
&&F^{\prime\prime}(t)\geq \ve B_0(t+R)^{-1}, \q t\geq 0,  \label{e3.15} \\
&&F^\prime(t)\geq \ve B_0\ln\lt(\f{t+R}{R}\rt), \q t\geq 0, \label{e3.16} \\
&&F(t)\geq C\ve B_0(t+R)\ln\lt(\f{t+R}{R}\rt),\q t\geq t_2, \label{3.14}
\eea
where $t_2=\max\{t_1,R(e^2-1)\}.$

Actually from \eqref{3.13}, \eqref{e3.15}, \eqref{e3.16} and \eqref{3.14}, we can deduce the blow up as in \cite{T03}. However for completion of our paper, we sketch the proof in the following. For simplicity, we set $R=1$.

Inserting \eqref{3.14} into \eqref{3.13}, one obtain the improvement for $F^{\prime\prime}(t)$
\be
F^{\prime\prime}(t)\geq C\ve^2B^2_0(t+1)^{-1}\ln(t+1), \q t\geq t_2.  \label{3.15}
\ee
Integrating \eqref{3.15} twice, we have
\be
F(t)\geq C\ve^2 B^2_0(t+1)\lt[\ln(t+1)\rt]^2 \q t\geq t_2.  \label{3.16}
\ee
Inserting \eqref{3.16} into \eqref{3.13}, we have
\be
F^{\prime\prime}(t)\geq C\ve^2 B^2_0(t+1)^{-2}\ln(t+1)F(t), \q t\geq t_2.  \label{3.17}
\ee
Multiplying both sides of \eqref{3.17} by $F^\prime(t)(\geq0),$ we have
\be
[(F^\prime(t))^2]^\prime\geq C\ve^2 B^2_0(t+1)^{-2}\ln(t+1)[(F(t))^2]^\prime. \label{3.18}
\ee
Integrating \eqref{3.18} from some $t_3\geq t_2$ to $t$, we have
\be
(F^\prime(t))^2\geq (F^\prime(t_3))^2+C \ve^2B^2_0\int^t_{t_3}(\tau+1)^{-2}\ln(\tau+1)[F^2(\tau)]^\prime d\tau. \nn
\ee
First choosing $t_3$ such that
\be
C\ve^2B^2_0\ln(t_3+1)= 1.  \label{e3.19}
 \ee
Then we have
\be
(F^\prime(t))^2\geq (F^\prime(t_3))^2+\f{1}{\ln(t_3+1)}\int^t_{t_3}(\tau+1)^{-2}\ln(\tau+1)[F^2(\tau)]^\prime d\tau. \nn
\ee
Using integration by parts, we have
\bea
(F^\prime(t))^2&\geq& (F^\prime(t_3))^2+\f{(t+1)^{-2}\ln(t+1)}{\ln(t_3+1)}F^2(t) \nn  \\
 &&-(t_3+1)^{-2}F^2(t_3)-\f{1}{\ln(t_3+1)}\int^t_{t_3}F^2(\tau)\f{1-2\ln(t+1)}{(t+1)^3} d\tau \nn \\
 &\geq&\f{(t+1)^{-2}\ln(t+1)}{\ln(t_3+1)}F^2(t) \nn \\
 &&+(F^\prime(t_3))^2-(t_3+1)^{-2}F^2(t_3). \label{3.19}
\eea
Here we have use the fact when $t\geq t_3,$ we have $1-2\ln(t+1)\leq 0.$
Noting that $F^\prime(t)$ is increasing and $F(0)=0,$ we have
\be
F(t_3)\leq t_3F^\prime(t_3).  \label{3.20}
\ee
Substituting \eqref{3.20} into \eqref{3.19}, we have
\bea
F^\prime(t)&\geq&  C(t+1)^{-1}(\ln(t+1))^\f{1}{2}F(t) \q t\geq t_3.  \nn
\eea
Integrating this from $t_3$ to $t$, one obtain
\be
\ln\f{F(t)}{F(t_3)}=\f{2}{3}C(\ln(t+1))^{3/2}-\f{2}{3}C(\ln(t_3+1))^{3/2}.  \nn
\ee
Choosing $t_4=2t^2_3$ and noting \eqref{3.16}, we have
\be
F(t)\geq C\ve^2B^2_0(t+1)^8,\q t\geq t_4.  \nn
\ee
Inserting this into \eqref{3.13}, we get
\be
F^{\prime\prime}(t)\geq C\ve B_0F(t)^{3/2}, \q t\geq t_4.  \nn
\ee
Integrating, as before, the above differential inequality, we get
\be
(F^\prime(t))^2\geq C\ve B_0\lt((F(t))^{5/2}-(F(t_4))^{5/2}\rt),  \q  t\geq t_4.  \nn
\ee
On the other hand, due to the nonnegative of $F^\prime(t)$ and $F^{\prime\prime}(t)$, we have
\be
F(t)\geq F(t_4)+F^\prime(t_4)(t-t_4)\geq F^\prime(t_4)(t-t_4)\geq F(t_4)\f{t-t_4}{t_4}.  \nn
\ee
Then choosing $t_5=3t_4$, we get
\be
F^\prime(t))\geq C\sqrt{\ve B_0}(F(t))^{5/4},\q t\geq t_5. \nn
\ee
If the lifespan $T_\ve \geq 2t_5$,  integrating the above inequality from $t_5$ to $T_\ve$ gives
\be
(F(t_5))^{-1/4}-(F(T_\ve))^{-1/4}\geq C\sqrt{\ve B_0}T_\ve. \label{3.25}
\ee
Noting that we have chosen $t_5=6t^2_3$ and the inequality \eqref{3.16} and \eqref{e3.19}, Then we have
\be
F(t_5)\geq C\ve^2B^2_0e^{\f{C}{\ve^2B^2_0}},  \nn
\ee
which combined with \eqref{3.25} can deduce that
\be
T_\ve\leq Ce^{\f{C}{\ve^2 B^2_0}}.  \nn
\ee
However if $T_\ve\leq 2 t_5$, then by our choice of $t_5$ and \eqref{e3.19}
\be
T_\ve\leq Ct^2_3\leq Ce^{\f{C}{\ve^2B^2_0}}.    \nn
\ee

For the general case, we need to adjust the function $G(r,t)$ in \eqref{e3.2}. Using Taylor's theorem, we have
\be
(\rho^\g-1)-\g(\rho-1)=\g(\g-1)\int^\rho_1\tau^{\g-2}(\rho-\tau)d\tau.  \nn
\ee
It is easy to see that
\be
\int^\rho_1\tau^{\g-2}(\rho-\tau)d\tau\geq C(\g)\varphi_\g(\rho),  \nn
\ee
where $C(\g)$ is a constant and $\varphi_\g$ is given by
\be
\varphi_\g(\rho)=\lt\{
\begin{aligned}
&(1-\rho)^\g,\q 0<\rho<\f{1}{2},  \nn  \\
&(\rho-1)^2,\q \f{1}{2}\leq \rho\leq 2, \nn  \\
&(\rho-1)^\g, \q \rho>2.  \nn
\end{aligned}\rt.
\ee
Then
\be
G(r,t)\geq C(\g)\int_{x>r}\varphi_\g(\rho)dx.  \nn
\ee

Young inequalities will be used in \eqref{3.12}. We still can get similar inequalities as \eqref{3.13}, \eqref{e3.15}, \eqref{e3.16} and \eqref{3.14} to prove the finite-time blow up, although the upper bound for the lifespan will be a little different. We omit the details.

This finishes the proof of Theorem1.2. \ef

\indent

{\bf Acknowledgement.} I want to express my heartfelt gratitude to my advisor Professor Huicheng Yin in Nanjing Normal University for his guidance about this subject. This work is partly proceeded when I am visiting Department of Mathematics in University of California, Riverside. So, I also want to express my thanks to my co-advisor, Professor Qi S. Zhang, in UCR for his constant encouragement.

\
\\

\end{document}